\documentclass[11pt]{amsart}
\usepackage{amssymb}
\usepackage{verbatim}

\begin{document}

\newtheorem{theorem}{Theorem}    
\newtheorem{proposition}[theorem]{Proposition}
\newtheorem{conjecture}[theorem]{Conjecture}
\def\theconjecture{\unskip}
\newtheorem{corollary}[theorem]{Corollary}
\newtheorem{lemma}[theorem]{Lemma}
\newtheorem{observation}[theorem]{Observation}
\theoremstyle{definition}
\newtheorem{definition}{Definition}
\newtheorem{remark}{Remark}
\newtheorem{question}{Question}
\def\thequestion{\unskip}
\newtheorem{example}{Example}
\def\theexample{\unskip}
\newtheorem{problem}{Problem}

\numberwithin{theorem}{section}
\numberwithin{definition}{section}
\numberwithin{equation}{section}
\numberwithin{remark}{section}

\def\jarrow{{\mathbf j}}
\def\Narrow{{\mathbf N}}
\def\Marrow{{\mathbf M}}

\def\reals{{\mathbb R}}
\def\torus{{\mathbb T}}
\def\integers{{\mathbb Z}}
\def\complex{{\mathbb C}\/}
\def\naturals{{\mathbb N}\/}
\def\distance{\operatorname{distance}\,}
\def\degree{\operatorname{degree}\,}
\def\kernel{\operatorname{kernel}\,}
\def\dim{\operatorname{dimension}\,}
\def\Span{\operatorname{span}\,}
\def\ZZ{ {\mathbb Z} }
\def\e{\varepsilon}
\def\p{\partial}
\def\rp{{ ^{-1} }}
\def\Re{\operatorname{Re\,} }
\def\Im{\operatorname{Im\,} }
\def\ov{\overline}
\def\bx{{\bf{x}}}
\def\eps{\varepsilon}
\def\lt{L^2}
\def\Tterm{T^\infty}
\def\Tnon{T^0}

\def\scriptx{{\mathcal X}}
\def\scriptj{{\mathcal J}}
\def\scriptr{{\mathcal R}}
\def\scripts{{\mathcal S}}
\def\scriptb{{\mathcal B}}
\def\scripta{{\mathcal A}}
\def\scriptk{{\mathcal K}}
\def\scriptd{{\mathcal D}}
\def\scriptp{{\mathcal P}}
\def\scriptl{{\mathcal L}}
\def\scriptv{{\mathcal V}}
\def\scripti{{\mathcal I}}
\def\scripth{{\mathcal H}}
\def\scriptm{{\mathcal M}}
\def\scripte{{\mathcal E}}
\def\scriptt{{\mathcal T}}
\def\scriptb{{\mathcal B}}
\def\scriptf{{\mathcal F}}
\def\scriptn{{\mathcal N}}
\def\frakg{{\mathfrak g}}
\def\frakG{{\mathfrak G}}

\author{Michael Christ}
\address{
        Michael Christ\\
        Department of Mathematics\\
        University of California \\
        Berkeley, CA 94720-3840, USA}
\email{mchrist@math.berkeley.edu}
\thanks{The author was supported by NSF grant 
DMS-040126}

\date{December 18, 2004. Revised February 15, 2005}

\title [Nonlinear Schr\"odinger equation]
{Power series solution \\ of a nonlinear Schr\"odinger equation}

\begin{abstract}
A slightly modified variant of the cubic periodic one-dimensional
nonlinear Schr\"odinger equation is shown to be well-posed,
in a relatively weak sense, in certain function spaces
wider than $L^2$.
Solutions are constructed as sums of infinite series of 
multilinear operators applied to initial data;
no fixed point argument or energy inequality are used.
\end{abstract}

\maketitle

\section{Introduction}
\subsection{The NLS Cauchy problem } 
The Cauchy problem for the one-dimensional periodic cubic
nonlinear Schr\"odinger equation is
\begin{equation} \label{nlsivp}
\left\{
\begin{aligned}
&iu_t + u_{xx} + \omega |u|^2 u=0
\\
&u(0,x)=u_0(x)
\end{aligned}
\right.
\tag{NLS}
\end{equation}
where $x\in\torus=\reals/2\pi\integers$,
$t\in\reals$, and the parameter $\omega$ equals $\pm 1$.
Bourgain \cite{bourgainperiodic} has shown this problem to be
wellposed in the Sobolev space $H^s$
for all $s\ge 0$, in the sense of uniformly continuous dependence on
the initial datum. 
In $H^0$ it is wellposed globally in time,
and as is typical in this subject, 
the uniqueness aspect of wellposedness is formulated in a certain auxiliary
space more restricted than $C^0([0,T], H^s(\torus))$, in which existence is
also established.
For $s<0$ it is illposed in the sense of uniformly continuous dependence 
\cite{burqgerardtzvetkov},
and is illposed in stronger senses \cite{cct3} as well.

The objectives of this paper are twofold. Firstly, we seek to establish the existence
of solutions for wider classes of initial data than $H^0$.
Secondly, we aim to develop an alternative method of solution.

The spaces of initial data considered here are the spaces $\scripth^{s,p}$
for $s\ge 0$ and $p\in [1,\infty]$,  defined as follows:
\begin{definition}
$\scripth^{s,p}(\torus)=\{f\in \scriptd(\torus): \langle \cdot\rangle^s
\widehat{f}(\cdot)\in\ell^p\}$.
\end{definition}
Here $\scriptd(\torus)$ is the usual space of distributions, and $\scripth^{s,p}$
is equipped with the norm
$\|f\|_{\scripth^{s,p}} = \|\widehat{f}\|_{\ell^{s,p}(\integers)}
= \big(\sum_{n\in\integers} \langle n\rangle^{ps}|\widehat{f}(n)|^p
\big)^{1/p}$.
We write $\scripth^p=\scripth^{0,p}$, and are mainly interested in these spaces since,
for $p>2$, they are larger function spaces than the borderline Sobolev space $H^0$ in which
\eqref{nlsivp} is already known to be wellposed.

\subsection{Motivations}
At least four concrete considerations motivate analysis of the Cauchy problem
in  these particular function spaces.
Firstly, 
$\scripth^p$ scales like $H^{s(p)}$ where $s(p)=-\tfrac12+\tfrac1p
\downarrow -\tfrac12$ as $p\uparrow \infty$, thus spanning the gap between
the optimal exponent $s=0$ for Sobolev space wellposedness,
and the scaling exponent $-\tfrac12$.

A second motivation is the work of
Kappeler and Topalov \cite{kappelertopalov1},\cite{kappelertopalov2},
who showed via an inverse scattering analysis that
the periodic KdV and mKdV equations
are wellposed for wider ranges of Sobolev spaces $H^s$ than had previously been known.
It is reasonable to seek a corresponding improvement for \eqref{nlsivp}, 
but this problem has been shown to be illposed in strong senses in $H^s$
for all $s<0$ \cite{cct3}.
Christ and Erdogan have investigated in unpublished work the 
relatively simple ``action variable''
portion of the inverse scattering theory relevant to \eqref{nlsivp},
and have found that for any
distribution in $\scripth^p(\torus)$ with small norm,
the sequence of gap lengths for the associated Dirac operator belongs
to $\ell^p$ and has comparable 
norm.\footnote{Having slightly better than bounded Fourier coefficients seems
to be a minimal condition for the applicability of this machinery,
since the eigenvalues for the free periodic Dirac system are  
equally spaced, and gap lengths for perturbations are to leading order
proportional to absolute values of Fourier coefficients of 
the perturbing potential.}
Thus 
$\scripth^p$ for $2<p<\infty$ may be a natural setting 
for inverse scattering theory for the Dirac operator relevant to
the periodic cubic nonlinear Schr\"odinger equation.

A related third motivation is the  goal 
of developing an alternative approach to the results of Kappeler and 
Topalov, independent of inverse scattering theory. NLS seems to be technically
simpler than mKdV or KdV, so it may be a reasonable starting point.
Fourthly, Gr\"unrock \cite{grunrock}
has proved wellposedness for the cubic nonlinear
Schr\"odinger equation in spaces analogous to $\scripth^{s,p}$, with $\torus$ 
replaced by $\reals$, and for other PDE in these function spaces, as well.


\subsection{Modified equation}
In order for the Cauchy problem to make any sense in $\scripth^p$
for $p>2$ it seems to be essential to modify the differential equation.
We consider
\begin{equation} \label{modnlsivp}
\left\{
\begin{aligned} 
&iu_t + u_{xx} + \omega \big( |u|^2-2\mu(|u|^2)) u=0
\\
&u(0,x)=u_0(x)
\end{aligned}
\right.
\tag{NLS$^*$}
\end{equation}
where 
\begin{equation}
\mu(|f|^2) = (2\pi)\rp\int_{\torus} |f(x)|^2\,dx
\end{equation}
equals the mean value of the absolute value squared of $f$.
In \eqref{modnlsivp}, 
$\mu(|u|^2)$ is shorthand for $\mu(|u(t,\cdot)|^2)=\|u(t,\cdot)\|_{\lt}^2$,
which is independent of $t$ for all sufficiently smooth solutions;
modifying the equation in this way merely introduces a unimodular scalar 
factor $e^{2i\mu t}$, where $\mu =\mu(|u_0|^2)$.
For parameters $p,s$ such that $\scripth^{s,p}$ is not embedded in $H^0$,
$\mu(|u_0|^2)$ is not defined for typical $u_0\in\scripth^{s,p}$, but 
of course the same goes for the function $|u_0(x)|^2$, and we will nonetheless
prove that the equation makes reasonable sense for such initial data.

The coefficient $2$ in front of $\mu(|u|^2)$ is
the unique one for which solutions depend continuously on initial data in
$\scripth^p$ for $p>2$.

\subsection{Conclusions}
Our main result is as follows. Recall that 
there exists a unique mapping $u_0\mapsto Su_0(t,x)$, 
defined for $u_0\in C^\infty$,
which for all sufficiently large $s$ extends to a uniformly continuous 
mapping from $H^s(\torus)$ to $C^0([0,\infty), H^s(\torus))\cap
C^1([0,\infty), H^{s-2}(\torus))$, such that $Su_0$ is a solution of the
modified Cauchy problem \eqref{modnlsivp}.
$C^\infty(\torus)$ is of course a dense subset of $\scripth^{s,p}$ for any $p\in[1,\infty)$.
\begin{theorem} \label{maintheorem}
For any $p\in[1,\infty)$, any $s\ge 0$,
and any $R<\infty$,
there exists $\tau>0$ for which
the solution mapping $S$ extends by continuity
to a uniformly continuous mapping from
the ball centered at $0$ of radius $R$ in 
$\scripth^{s,p}(\torus)$ 
to $C^0([0,\tau],\scripth^{s,p}(\torus))$.
\end{theorem}

For the unmodified equation this has the following consequence.
Denote by $H^0_c=H^0_c(\torus)$ the set of all $f\in H^0$ such
that $\|f\|_{\lt}=c$.
Denote by $S'u_0$ the usual solution \cite{bourgainperiodic} of the unmodified Cauchy problem
\eqref{nlsivp} with initial datum $u_0$, for $u_0\in H^0$.
\begin{corollary}
Let $p\in[1,\infty)$ and $s\ge 0$.
For any $R<\infty$ there exists $ \tau>0$ such that
for any finite constant $c>0$,
the mapping $H^0_c\owns u_0\mapsto S'u_0$ 
is uniformly continuous as a mapping from $H^0_c$ intersected with the ball
centered at $0$ of radius $R$ in $\scripth^{s,p}$, 
equipped with the $\scripth^{s,p}$ norm,
to $C^0([0,\tau],\scripth^{s,p}(\torus))$.
\end{corollary}

The unpublished result of the author and Erdogan
shows that for initial data in $\lt$, for which the solution is known to exist
globally in time, $\|u(t)\|_{\scripth^p}\le C\|u_0\|_{\scripth^p}$
uniformly for all $t\in[0,\infty)$,
provided that $\|u_0\|_{\scripth^p}$ is sufficiently small.
This result, once published, will combine with Theorem~\ref{maintheorem}
to yield global wellposedness for sufficiently small data.

The following result quantifies the relation between the nonlinear evolution
\eqref{modnlsivp} ond the corresponding linear Cauchy problem
\begin{equation} \label{linearevolution}
\left\{
\begin{aligned}
&iv_t + v_{xx}=0
\\
&
v(0,x)=u_0(x).
\end{aligned}
\right.
\end{equation}
\begin{proposition} \label{prop:weaksmoothing}
Let $R<\infty$ and $p\in[1,\infty)$.
Let $q>p/3$ also satisfy $q\ge 1$.
Then there exist $\tau,\eps>0$  and $C<\infty$ such that
for any initial datum $u_0$ satisfying $\|u_0\|_{\scripth^p}\le R$,
the solutions $u=Su_0$ of \eqref{modnlsivp} and $v$ of \eqref{linearevolution}
satisfy
\begin{equation}
\|u(t,\cdot)-v(t,\cdot)\|_{\scripth^q}
\le Ct^\eps\ \ \text{for all } t\in[0,\tau].
\end{equation}
\end{proposition}
\noindent
Here $u$ the solution defined by approximating $u_0$ by elements of $C^\infty$
and passing to the limit. Thus for $p>1$ the nonlinear terms are in a sense
smoother than the linear evolution.

Our next result 
indicates that the function
$u(t,x)$ defined by the limiting procedure of Theorem~\ref{maintheorem}
is a solution of the differential equation in a more natural sense
than merely being a limit of smooth solutions.
Define Fourier truncation operators
$T_N$, acting on $\scripth^{s,p}(\torus)$,
by $\widehat{T_Nf}(n)=0$ for all $|n|>N$, and $=\widehat{f}(n)$
whenever $|n|\le N$. 
$T_N$ acts also on functions $v(t,x)$ by acting on $v(t,\cdot)$
for each time $t$ separately.
We denote by $S(u_0)$ the limiting function 
whose existence, for nonsmooth $u_0$,
is established by Theorem~\ref{maintheorem}.
\begin{proposition} \label{prop:nonlineardistribution}
Let $p\in[1,\infty)$, $s\ge 0$, and $u_0\in\scripth^{s,p}$. 
Write $u=S(u_0)$. 
Then for any $R<\infty$ there exists $\tau>0$ such that whenever
$\|u_0\|_{\scripth^{s,p}}\le R$,
$\scriptn u(t,x) = (|u|^2-2\mu(|u|^2))u$ exists in the sense that
\begin{equation} \label{Nlimit}
\lim_{N\to\infty} \scriptn(T_N u)(t,x)
\text{ exists in the sense of distributions in }
C^0((0,\tau),\scriptd'(\torus)).
\end{equation}
Moreover if $\scriptn(u)$ is interpreted as this limit, then
$u=S(u_0)$ satisfies \eqref{modnlsivp} in the sense of distributions
in $(0,\tau)\times\torus$.
\end{proposition}
\noindent More generally, the same holds for any sequence of 
Fourier multipliers of the form $\widehat{T_\nu f}(n) = m_\nu(n)\widehat{f}(n)$
where each sequence $m_\nu$ is finitely supported,
$\sup_\nu\|m_\nu\|_{\ell^\infty} <\infty$,
and $m_\nu(n)\to 1$ as $\nu\to\infty$ for each $n\in\integers$;
the limit is of course independent of the sequence $(m_\nu)$.
Making sense of the nonlinearity via this limiting procedure 
is connected with general theories of multiplication of
distributions \cite{biagioni},\cite{colombeau},
but the existence here of the limit over all sequences $(m_\nu)$
gives $u$ stronger claim to the title of solution than 
in the general theory.

Unlike the fixed point method, our proof
yields no uniqueness statement corresponding to these existence results.
But this failing is unavoidable;
for all $p>2$, solutions of the Cauchy problem in the class 
$C^0([0,\tau],\scripth^p)$,
in the sense of Proposition~\ref{prop:nonlineardistribution},
are not unique \cite{nonuniqueness}. 

\subsection{Method}
Define the partial Fourier transform
\begin{equation}
\widehat{u}(t,n) = (2\pi)\rp\int_{\torus} e^{-inx}u(t,x)\,dx.
\end{equation}
Our approach is to regard the partial differential equation as an infinite
coupled nonlinear system of ordinary differential equations for these
Fourier coefficients, to express the solution as a power series in the
initial datum
\begin{equation}
\widehat{u}(t,n)
= \sum_{k=0}^\infty \hat A_k(t)(\widehat{u_0},\cdots,\widehat{u_0})
\end{equation}
where each $\hat A_k(t)$ is a bounded multilinear operator\footnote{
Throughout the discussion
we allow multilinear operators to be either conjugate linear or linear 
in each of their arguments, independently.}
from a product of $k$
copies of $\scripth^{s,p}$ to $\scripth^{s,p}$, to show that the individual terms
$\hat A_k(t)(\widehat{u_0},\cdots,\widehat{u_0})$ 
are well-defined, and to show that the formal
series converges absolutely in $C^0(\reals,\scripth^{s,p})$ to a solution
in the sense of \eqref{Nlimit}.
The case $s\ge 0$ follows from a very small modification of the analysis for $s=0$,
so we discuss primarily $s=0$, indicating the necessary modifications for $s>0$
at the end of the paper.

The analysis is quite elementary,
much of the paper being devoted to setting up the definitions and notation
required to describe the operators $\hat A_k(t)$. 
A single number theoretic fact enters the discussion: the number of factorizations
of an integer $n$ as a product of two integer factors is $O(n^\delta)$ as $n\to\infty$,
for all $\delta>0$; this same fact was used in a more sophisticated way 
by Bourgain \cite{bourgainperiodic}.

\medskip
The author is grateful to J.~Bourgain, C.~Kenig, H.~Koch, and D.~Tataru
for invitations to conferences that stimulated this work,
and to Betsy Stovall for thorough proofreading of the manuscript.

\section{A system of coupled ordinary differential equations} \label{section:ode}
\subsection{General discussion}
Define
\begin{equation}
\sigma(j,k,l,n) = n^2-j^2+k^2-l^2.
\end{equation}
It factors as
\begin{equation}
\sigma(j,k,l,n) = 2(n-j)(n-l) = 2(k-l)(k-j)
\ \text{provided that $j-k+l=n$}.
\end{equation}

Written in terms of Fourier coefficients $\widehat{u}_n(t) = \widehat{u}(t,n)$, 
the equation $iu_t + u_{xx} + \omega \big(|u|^2-2\mu(|u|^2) \big)u=0$
becomes
\begin{equation} \label{firstode}
i\frac{d \widehat{u}_n}{dt} -n^2 \widehat{u}_n
+ \omega \sum_{j-k+l=n} \widehat{u}_j\overline{\widehat{u}_k}\widehat{u}_l
-2\omega \sum_m |\widehat{u}_m|^2\widehat{u}_n
=0.
\end{equation}
Here the first summation is taken over all $(j,k,l)\in\integers^3$ 
satisfying the indicated identity, and the second over all $m\in\integers$.
Substituting
\begin{equation} 
a_n(t) = e^{in^2 t}\widehat{u}(t,n),
\end{equation}
\eqref{firstode} becomes
\begin{equation} \label{secondode}
\frac{d a_n}{dt} = i\omega \sum_{j-k+l=n}^* a_j \bar a_k a_l e^{i \sigma(j,k,l,n)t}
-i\omega |a_n|^2 a_n.
\end{equation}
where the notation $\sum_{j-k+l=n}^*$ means that the sum is taken over all
$(j,k,l)\in\integers^3$ for which neither $j=n$ nor $l=n$. 
This notational convention will be used throughout the discussion.
The effect of the term $-2\omega\mu(|u|^2)u$ in the modified differential
equation \eqref{modnlsivp} is to cancel out a term
$2i\omega(\sum_m |a_m|^2)a_n$, which would otherwise appear on the right-hand side
of \eqref{secondode}. 

Reformulated as an integral equation, \eqref{secondode} becomes
\begin{equation} \label{firstintegraleqn}
a_n(t) = a_n(0)
+ i\omega\sum_{j-k+l=n}^* 
\int_0^t a_j(s) \bar a_k(s) a_l(s)e^{i\sigma(j,k,l,n)s}\,ds
-i\omega\int_0^t |a_n(s)|^2 a_n(s)\,ds.
\end{equation}
However, in deriving \eqref{firstintegraleqn}
from \eqref{secondode}, we have interchanged the
integral over $[0,t]$ with the summation over $j,k,l$ without any justification.
\eqref{firstintegraleqn} is fully equivalent to
\begin{equation} \label{ufirstintegraleqn}
\widehat{u}(t,n)
= \widehat{u_0}(n)
-in^2 \int_0^t \widehat{u}(s,n)\,ds
+ i\omega\sum_{j-k+l=n}^*
\int_0^t 
\widehat{u}(s,j)
\overline{\widehat{u}(s,k)}
\widehat{u}(s,l)
\,ds
-i\omega\int_0^t |\widehat{u}(s,n)|^2\widehat{u}(s,n)\,ds.
\end{equation}

Substituting for $a_j,a_k,a_l$ in the right-hand side of \eqref{firstintegraleqn}
by means of the equation itself yields
\begin{align}
a_n(t) & = a_n(0)
+ i\omega\sum_{j-k+l=n}^* 
a_j(0)\bar a_k(0)\bar a_l(0)\int_0^t e^{i\sigma(j,k,l,n)s}\,ds
-i\omega |a_n(0)|^2a_n(0)\int_0^t 1\,ds
\label{anform}
\\
& \qquad\qquad+ \text{ additional terms}.
\notag
\\
&= a_n(0)\big(1-i\omega t|a_n(0)|^2 \big) 
+\tfrac12\omega \sum_{j-k+l=n}^* \frac{a_j(0)\bar a_k(0) a_l(0)}
{(n-j)(n-l)} (e^{i(n^2-j^2+k^2-l^2)t}-1)
\notag
\\
&\qquad\qquad+\text{ additional terms.}
\notag
\end{align}
We recognize $1-i\omega t|a_n(0)|^2$ as a Taylor polynomial for
$\exp(-i|a_n(0)|^2 t)$, but for our purposes it will not be necessary to
exploit this by recombining terms, 
and in particular we will not exploit the coefficient $i$ which makes
this exponential unimodular.

\subsection{A sample term}
One representative additional term is
\begin{multline} \label{presample}
(i\omega)^4 
\sum_{j_1-j_2+j_3=n}^*
\,\,
\sum_{m^1_1-m^1_2+m^1_3=j_1}^*
\,\,
\sum_{m^2_1-m^2_2+m^2_3=j_2}^*
\,\,
\sum_{m^3_1-m^3_2+m^3_3=j_3}^*
\\
\int_{0\le r_1,r_2,r_3\le s\le t}
a_{m^1_1}(r_1)\bar a_{m^1_2}(r_1) a_{m^1_3}(r_1)
\bar a_{m^2_1}(r_2) a_{m^2_2}(r_2) \bar a_{m^2_3}(r_2)
a_{m^3_1}(r_3)\bar a_{m^3_2}(r_3) a_{m^3_3}(r_3)
\\
e^{i\sigma(j_1,j_2,j_3,n)s}
e^{i\sigma(m^1_1,m^1_2,m^1_3,j_1)r_1}
e^{-i\sigma(m^2_1,m^2_2,m^2_3,j_2)r_2}
e^{i\sigma(m^3_1,m^3_2,m^3_3,j_3)r_3}
\,dr_1\,dr_2\,dr_2\,ds.
\end{multline}
Substituting via \eqref{firstintegraleqn} for each coefficient $a$ 
yields a main term
\begin{equation} \label{firstsample}
(i\omega)^4 
\sum_{(m^i_k)_{1\le i,k\le 3}}^*
\scripti( t, (m^i_k)_{1\le i,k\le 3} )
\prod_{i,j=1}^3 a^*_{m^i_j}(0)
\end{equation}
plus higher-degree terms,
where the superscript $*$ indicates here that the sum is taken over only
certain $(m^i_k)_{1\le i,k\le 3}\subset\integers^9$ (more precisely,
over most of a copy of $\integers^8$ affinely embedded in $\integers^9$), 
$a^*_{m^i_j}(0) = a_{m^i_j}(0)$ if $i+j$ is even
and $=\overline{a_{m^i_j}(0)}$ if $i+j$ is odd,
and
\begin{equation}
\scripti(t,(m^i_k)_{1\le i,k\le 3})
=
\int\limits_{0\le r_1,r_2,r_3\le s\le t}
e^{i\theta(t,s,r_1,r_2,r_3,\{m^i_j: 1\le i,j\le 3\})}
\,dr_1\,dr_2\,dr_2\,ds,
\end{equation}
with 
\begin{equation}
\theta( t,s,r_1,r_2,r_3,(m^i_j)_{ 1\le i,j\le 3} )
=
\sigma(j_1,j_2,j_3,n)s
+ \sum_{i=1}^3 (-1)^{i+1}
\sigma(m^i_1,m^i_2,m^i_3,j_i)r_i;
\end{equation}
here $j_1,j_2,j_3,n$ are defined as functions of $(m^i_j)$ 
by the equations governing the sums in \eqref{presample}.
Continuing in this way yields formally
an infinite expansion for the sequence $(a_n(t))_{n\in\integers}$
in terms of multilinear expressions in the initial datum $(a_n(0))$.
This expansion is doubly infinite; the single (and relatively simple)
term \eqref{firstsample}
is for instance an infinite sum over most of a copy of $\integers^8$ for each $n$.

The discussion up to this point has been purely formal, 
with no justification of convergence.
In the next section we will begin to describe the terms in this expansion
systematically.

\section{Trees} \label{section:trees}

On a formal level $a(t) = (a_n(t))_{n\in\integers}$
equals an infinite sum $\sum_{k=1}^\infty A_k(t)(a(0),a(0),a(0),\cdots)$
where each $A_k(t)$ is a sum of finitely many multilinear operators, each of degree $k$.
We now describe a class of trees which will be used both to name,
and to analyze, these multilinear operators.

\begin{definition}
A tree $T$ is a finite partially ordered set with the following properties:
\begin{enumerate}
\item
Whenever $v_1,v_2,v_3,v_4\in T$ and $v_4\le v_2\le v_1$ and $v_4\le v_3\le v_1$,
then either $v_2\le v_3$ or $v_3\le v_2$.
\item
There exists a unique element $v_0\in T$ satisfying $v_0\ge v$ for all $v\in T$.
\item
Each $v\in T$ has either three children, or no children; 
$w$ is said to be a child of $v$
if $w<v$ and if there exists no $u\in T$ satisfying $w<u<v$.
\item For each $v \in T$ there is given an element of $\{\pm 1\}$,
denoted $\pm_v$.
\end{enumerate}
\end{definition}
 
\begin{definition}
Elements of $T$ are called nodes. A terminal node is one with zero children.
The maximal element of $T$ is called its root node.
For any $u\in T$, $T_u=\{v\in T: v\le u\}$ is a tree, with root node $u$.
$\Tterm$ denotes the set of all terminal nodes of $T$,
while $\Tnon=T\setminus\Tterm$ denotes the set $\Tnon$ of all non-terminal nodes.
The three children of any $v\in\Tnon$ are denoted by $(v,1),(v,2),(v,3)$.
\end{definition}

The number $|T|$ of nodes of a tree is of the form $1+3k$ for some nonnegative
integer $k$. The number of terminal nodes is then 
\begin{equation}
|\Tterm| = 1+2k =  \tfrac23 |T| + \tfrac13.
\end{equation}

\begin{definition}
An ornamented tree is a tree $T$ together with the following additional structure:
\begin{enumerate}
\item
Associated to each node $v\in T$ is copy of $\integers$, indexed
by the variable $j_v$.
\item
There is given a partition of the set of all non-terminal nodes of $T$ into two
disjoint classes, called simple nodes and general nodes. Terminal nodes
are neither simple nor general.
\item
For each non-terminal node $v\in \Tnon$, 
and each $i\in\{1,2,3\}$ such that the child $(v,i)$ is non-terminal,
there is given a coefficient $\eps_{v,i}\in\{-1,0,1\}$.
\item
Associated to each node $v\in T$ is a $\integers$-valued function $\rho_v$
of $\jarrow=(j_u)_{ u\in T}$, defined by
\begin{equation}
\rho_v(\jarrow) = 0 \text{ if $v\in \Tterm$}
\end{equation}
and
\begin{equation}
\rho_v(\jarrow)
= \sigma(j_{(v,1)},j_{(v,2)},j_{(v,3)},j_v)
+ \sum_{i=1}^3 \eps_{v,i}\rho_{(v,i)}
\text{ if $v\in \Tnon$}.
\end{equation}
\end{enumerate}
\end{definition}
$\rho_v(\jarrow)$ actually depends only on $\{j_u,\eps_{u,i}: u\le v\}$.
We will use the symbol $T$ to denote both the ornamented tree and the underlying tree,
and will often write $\rho_v$ instead of $\rho_v(\jarrow)$.

\begin{definition}
Let $T$ be a tree.
$\scriptj(T)\subset\integers^T$
denotes the set of all
$\jarrow=(j_v)_{ v\in T}$ satisfying the restrictions
\begin{align}
&j_v = j_{(v,1)}-j_{(v,2)} + j_{(v,3)}
 \text{ for every $v\in \Tnon$}
\label{cyclicity}
\\
&\{j_v,j_{(v,2)}\}\cap\{j_{(v,1)},j_{(v,3)}\}=\emptyset
\ \text{ for every general node $v\in \Tnon$}
\label{exclusion}
\\
&j_v = j_{(v,i)} \ \text{ for all } i\in\{1,2,3\}
\ \text{ for every simple node $v\in \Tnon$}.
\label{simplicity}
\end{align}
\end{definition}

\begin{definition}
Let $T$ be any tree.
$\sigma_w:\scriptj(T)\to\integers$ denotes the function 
$\sigma_w(\jarrow)=0$
if $w$ is terminal, 
and $\sigma_w(\jarrow) = j_w^2-j_{(w,1)}^2 + j_{(w,2)}^2-j_{(w,3)}^2$
if $w$ is non-terminal.
\end{definition}

Let $\delta,c_0>0$
be sufficiently small positive numbers, to be chosen later.
The following key definition involves these quantities.

\begin{definition}
Given an ornamented tree $T$ and $\jarrow\in\scriptj(T)$, we say that
a node $v\in T$ is frozen if $v$ is non-terminal and
\begin{equation} \label{frozennode}
|\rho_v(\jarrow)|
\le c_0|\sigma(j_{(v,1)},j_{(v,2)},j_{(v,3)},j_v)|^{1-\delta}.
\end{equation}
If $v$ is not frozen, then $v$ is said to be alive.
A non-terminal node $v$ is said to be exceptional if $\rho_v(\jarrow)=0$. 
\end{definition}
Whether $v$ is frozen depends on the values of $j_u$ for all nodes $u\le v$,
as well as on $\eps_{u,i}$ for all non-terminal $u\le v$,
not merely on the structure of $T$; a non-terminal node will be frozen
for some $\jarrow$, but alive for others. 
Thus it would be more felicitous to say that a pair $(v,\jarrow)$
is frozen, rather than a node $v$.

Exceptional nodes are of course frozen.
If $v\in\Tnon$ is a general node
all three of whose children of $v$ are terminal, then $v$ cannot be exceptional, for
$\rho_v = \sigma(j_{(v,1)},j_{(v,2)},j_{(v,3)},j_v)
=2(j_v-j_{(v,1)})(j_v-j_{(v,3)})$ cannot vanish, by \eqref{exclusion}.
But if $v$ has at least one non-terminal child,
then nothing prevents $\rho_v$ from vanishing, and if $v$ is a simple node
all of whose children are terminal, then $v$ is certainly exceptional.

\begin{definition}
A weathered ornamented tree $(T,T')$
is an ornamented tree $T$ together with a subset $T'\subset \Tnon$
and the collection 
\begin{equation}
\scriptj(T,T') =\{\jarrow\in\scriptj(T): 
\text{ $v\in T$ is frozen if and only if $v\in T'$.}\}
\end{equation}
\end{definition}

\section{Multilinear operators associated to trees} \label{section:treeops}

\begin{definition}
Let $T$ be any tree, and let $t\in\reals$. 
The associated tree coefficients are
\begin{equation}
\scripti_T(t,\jarrow)
= \int_{\scriptr(T,t)}
\prod_{u\in \Tnon} e^{\pm_u i\omega\sigma_u(\jarrow) t_u}\,dt_u
\end{equation}
where
\begin{equation}
\scriptr(T,t)=\{(t_u)_{u\in\Tnon}:
0\le t_u\le t_{u'}\le t \text{ whenever } u\le u'\}.
\end{equation}
\end{definition}
\noindent


The following upper bounds for the coefficients $\scripti_T(t,\jarrow)$
are the only information concerning them that will be used in the analysis. 
\begin{lemma} \label{lemma:treecoefficientbound}
Let $T$ be any tree, and let 
$\jarrow\in\scriptj(T)$.
Then for all $t\in[0,1]$,
\begin{align}
& |\scripti_T(t,\jarrow)| \le t^{|\Tnon|}
\\
\intertext{and}
&|\scripti_T(t,\jarrow)| \le
2^{|T|}
\sum_{(\eps_{u,i})}
\prod_{w\in  \Tnon} \langle\rho_w(\jarrow)\rangle^{-1}.
\label{maincoefficientbound}
\end{align}
\end{lemma}
The notation $\langle x\rangle$ means  $(1+|x|^2)^{1/2}$. 
The sum here is taken over all of the $3^{|\Tnon|}$
possible choices of $\eps_{u,i}\in\{0,1,-1\}$;
these choices in turn determine the quantities $\rho_w$.
This lemma will be proved in \S\ref{section:treecoeffs}.

\begin{definition}
Let $T$ be any ornamented tree. The tree operator $\scripts_T(t)$ associated to 
$T$ is for each $t\in\reals$ the multilinear operator
that maps $(x_v)_{v\in\Tterm}$, where each $x_v$ is a sequence of complex 
numbers, to the sequence of complex numbers
\begin{equation}
\scripts_T(t)(x_v)_{ v\in \Tterm}(n) = \sum_{\jarrow\in\scriptj(T): j_{v_0}=n}
\scripti_T(t,\jarrow)
\prod_{w\in \Tterm}x_w(j_w).
\end{equation}
\end{definition}

\section{Formalities} 

With all these definitions and notations in place, we can finally 
formulate the conclusion of the discussion in \S\ref{section:ode};
proofs will be supplied later.
\begin{proposition} \label{prop:formalexpansion}
The recursive procedure indicated in \S\ref{section:ode}
yields a formal expansion 
\begin{equation} \label{partofformalsolutionseries}
a(t) = \sum_{k=1}^\infty A_k(t)(a(0),a(0),\cdots),
\end{equation}
where each $A_k(t)$ is a multilinear operator of the form
\begin{equation} \label{formalsolutionseries}
A_k(t) = \sum_{|T|=3k+1} c_T \scripts_T(t)
\end{equation}
where the scalars $c_T\in\complex$ satisfy
$|c_T|\le C^{k}$ for some finite constant $C$. 
The sum in \eqref{formalsolutionseries}
is taken over all ornamented trees $T$ of the indicated cardinalities.
There exists a finite positive constant $c_0$ such that whenever
$a(0)\in\ell^1$, the multiply infinite series $\sum_k A_k(t)(a(0),\cdots)$
converges absolutely to a function in $C^0([0,\tau],\ell^1)$
provided that $\tau\|a(0)\|_{\ell^1}\le c_0$.
\end{proposition}

By this last statement we mean that
$\sum_{\jarrow\in\scriptj(T)} |\scripti_T(t,\jarrow)|\prod_{w\in\Tterm}|a(0)(j_w)|$
converges absolutely for each ornamented tree $T$,
and that if its sum is denoted by $\scripts_T^*(a(0),a(0),\cdots)(t)$
then the resulting series
$\sum_{k=1}^\infty \sum_{|T|=3k+1} |c_T\scripts_T^*(a(0),a(0),\cdots)(t)|$
likewise converges.S
All but the last sentence follows from the discussion in \S\ref{section:ode}
and the definitions in \S\S\ref{section:trees},{section:treeops}.

The operators $\scripts_T$ and coefficients $c_T$
were defined so that the following holds automatically.
\begin{lemma} \label{automaticsolution}
There exists $c>0$ with the following property.
Let $\widehat{u_0}$ be any numerical sequence 
and define $a(0)(n) = \widehat{u_0}(n)$.
Suppose that the infinite series
defining $\scripts_T^*(a(0),a(0),\cdots)(t)$ converges absolutely and uniformly
for all $t\in[0,\tau]$ and that its sum is $O(c^{|T|})$,
uniformly for every ornamented tree $T$.
Define $a(t)$ to be the sequence $\sum_{k=1}^\infty A_k(t)(a(0),a(0),\cdots)$.
Then $a$ satisfies the integral equation 
\eqref{firstintegraleqn} for $t\in[0,\tau]$. 
Moreover the function $u(t,x)$ defined by
$\widehat{u}(t,n) = e^{-in^2 t}a(t,n)$
is a solution of the modified Cauchy problem \eqref{modnlsivp} in the corresponding
sense \eqref{ufirstintegraleqn}.
\end{lemma}

The main estimate in our analysis is as follows.
\begin{proposition} \label{prop:treeoperatorbound}
Let $p\in (1,\infty)$. Then for any exponent $q>\frac{p}{|\Tterm|}$
satisfying also $q\ge 1$,
there exist $\eps>0$ and $C<\infty$
such that for all ornamented trees $T$
and all sequences $x_v\in\ell^1$,
\begin{equation} \label{eq:treeoperatorbound}
\|\scripts_T(t)(x_v)_{v\in\Tterm}\|_{\ell^q}
\le (Ct^\eps)^{|\Tterm|} \prod_{v\in \Tterm}\|x_v\|_{\ell^p}.
\end{equation}
\end{proposition} 

Propositions~\ref{prop:formalexpansion}
and \ref{prop:treeoperatorbound}  and Lemma~\ref{automaticsolution}
will be proved in subsequent sections.
Together, they give:
\begin{corollary}
Let $p\in[1,\infty)$.
For any $R<\infty$
there exists $\tau>0$
such that the solution mapping $u_0\mapsto u(t,\cdot)$
for the modified Cauchy problem \eqref{modnlsivp},
initially defined for all sufficiently smooth $u_0$,
extends by uniform continuity to a real analytic mapping 
from $\{u_0\in \scripth^p: \|u_0\|_{\scripth^p}\le R\}$
to $C^0([0,\tau],\scripth^p(\torus))$. 
\end{corollary}

We emphasize that analytic dependence on $t$ is not asserted.

\section{Tree coefficient bound} \label{section:treecoeffs}

\begin{proof}[Proof of Lemma~\ref{lemma:treecoefficientbound}]
The first conclusion of the lemma
holds simply because $|\scripti_T(t,\jarrow)| \le |\scriptr(T,t)|$.
The proof of the main conclusion \eqref{maincoefficientbound}
proceeds recursively. In step 1 we integrate with respect to
$t_v$ for certain nodes $v$, holding fixed all other coordinates $t_w$ in the integral
defining $\scriptr(T,t)$.
Specifically, we
hold fixed the coordinate $t_v$ whenever at least one child of $v$ is not terminal.
We also fix $t_v$ for every simple node $v$ having only terminal children.
The former coordinates $t_v$, and underlying nodes $v$, are said to be temporarily fixed; 
the latter coordinates and nodes are said to be permanently fixed.
No other coordinates are fixed at this step.

When $|T|=1$ there is nothing to prove. 
Otherwise there must always exist at least one node, all of whose children are terminal
If there exists such a node which is also general, then
at least coordinate $t_v$ is not fixed.
The subset, or slice, of $\scriptr(T,t)$ defined by setting each of the fixed coordinates
equal to some constant is either empty, or takes the product form
$\times_{u \text{ not fixed}} [0,t_{u^*}]$,
where $u^*$ denotes the parent of $u$.
Integrating over this slice with respect to all of the non-fixed coordinates thus yields
\begin{equation*}
\prod_{w} e^{\pm_w i \omega \sigma_w t_w}
\prod_u \int_0^{t_{u^*}} e^{\pm_u i\omega\sigma_u t_u}\,dt_u,
\end{equation*}
where the first product is taken over all fixed $w\in \Tnon$, 
and the second over all remaining non-fixed $u\in \Tnon$.

None of the quantities $\sigma_u$ can vanish in step 1,
since a general node having only terminal children can never be exceptional,
by \eqref{exclusion}.
Therefore the preceding expression equals
\begin{equation*}
\prod_{w} e^{\pm_w i \omega \sigma_w t_w}
\prod_u 
(\pm_u i\omega\sigma_u)^{-1} \big(e^{\pm i \omega \sigma_u t_{u^*}} -1\big).
\end{equation*}
This may be expanded as a sum of $2^N$ terms, where $N$ is the number of
non-fixed nodes in $\Tnon$.
Each of these terms has the form
\begin{equation*}
\pm\prod_{w} e^{\pm_w i \omega \sigma_w t_w}
\prod_u 
(i\omega\sigma_u)^{-1} e^{\eps_u i \omega \sigma_u t_{u^*}} 
\end{equation*}
for some numbers $\eps_u\in\{0,1,-1\}$.

The other possibility in step 1 is that $|T|>1$, but every nonterminal node that has
only terminal children is simple.
In that case all coordinates $t_v$ are fixed at step 1, no integration
is performed, and we move on to step 2.

We now carry out step 2. If a node $v$ 
was permanently fixed at step 1 then it remains fixed for all subsequent steps;
we never integrate with respect to $t_v$. More generally, any node that is
permanently fixed at any step remains fixed through all subsequent steps.
On the other hand, once we've integrated with respect to some $t_w$, then the
node $w$ is also removed from further consideration.
All other nodes remain active, including those temporarily fixed in step 1. 
Denote by $T_1$ the set of all nodes that are active after the completion of step 1.

$T_1$ is itself a tree.
There is an associated subset $\scriptr_{T_1}$
of $\{(t_w: w\in T_1)\}$, defined by the inequalities $0\le t_w\le t_{w'}\le t$
whenever $w\le w'$, and also by $t_u\le t_w$ if $u\le w$ and $u$ was 
permanently fixed in step 1.
To each node $w\in T_1$ is associated a modified phase $\sigma^{(2)}_w$,
defined to be $\sigma_w + \sum_{i} \eps_{(w,i)} \sigma_{(w,i)}$,
where the sum is taken over all $i\in\{1,2,3\}$ such that we
integrated with respect to $t_{(w,i)}$ in the first step.

A node $w$ is permanently fixed at the second step if $w$ is terminal
in $T_1$ and satisfies $\sigma^{(2)}_w=0$.
A node $w\in T_{1}$ is temporarily fixed at the second step if $w$ is
not terminal in $T_{1}$.
We now integrate $\prod_{w\in T_1} e^{\pm i\omega\sigma^{(2)}_w(t_w)}$
over $\scriptr_{T_1}$
with respect to $t_u$ for all $u\in T_1$ that are neither temporarily 
nor permanently fixed. 
As in step 1, this integral has a product structure,
and $2^{N_2}$ terms are obtained, where $N_2$ is the number of 
variables with respect to which we integrate.

In step 3 we consider the tree $T_2$ consisting of all $w\in T_1$ that
were temporarily fixed in step 2. Associated to $T_2$
is a set $\scriptr_{T_2}$, and associated to each node $v\in T_2$
is a modified phase $\sigma^{(3)}_w =
\sigma^{(2)}_w + \sum_i\eps_{(w,i)}\sigma^{(2)}_{(w,i)}$,
the sum being taken over all $i\in\{1,2,3\}$ such that
$(w,i)$ was not fixed in step 2.
A node $v\in T_2$ is then permanently fixed if it is terminal in $T_2$
and $\sigma^{(3)}_v=0$.
$v\in T_2$ is temporarily fixed if it is not terminal in $T_2$.
We then integrate with respect to $t_v$ for all $v\in T_2$
that are neither temporarily nor permanently fixed. 

This procedure terminates after  finitely many steps, when for each node 
$v\in \Tnon$,
either $v$ has become permanently fixed, or we have integrated with respect to $t_v$.
This yields a sum of at most 
$2^{|\Tnon|}$ terms. Each term arises from some particular 
choice of the parameters $\eps_{u,i}$, and is expressed as an integral
with respect to $t_v$ for all nodes $v\in \Tnon$ that were
permanently fixed at some step; the vector $(t_v)$ indexed by all such $v$
varies over a subset of $[0,t]^M$ where $M$ is the number of such $v$.
At step $n$, each integration with respect to some $t_u$ yields
a factor of $(\sigma^{(n)}_u)^{-1}$, multiplied by some unimodular factor;
recall that $\sigma^{(n)}_u\ne 0$, since otherwise $u$ would have
been permanently fixed.

Thus for each term we obtain an upper bound of $\prod_u|\rho_u|^{-1}$,
where the product is taken over all nonexceptional nodes $u$;
this bound must still be 
integrated with respect to all $t_w$ where $w$ ranges over all the 
exceptional nodes. Each such coordinate $t_w$ is restricted to $[0,t]$.
Thus we obtain a total bound
\begin{equation}
|\scripti(t,\jarrow)|\le \sum_{(\eps_{u,i})}
t^M \prod_{w\in \Tnon}^* |\rho_w(\jarrow)|^{-1}
\end{equation}
where for each  $(\eps_{u,i})$, $M=M((\eps_{u,i}))$ 
is the number of exceptional nodes encountered
in this procedure, that is, the number of permanently fixed nodes, 
and where for each $(\eps_{u,i})$, $\prod_{w\in \Tnon}^*$
denotes the product over all nonexceptional nodes $w\in\Tnon$
that are nonexceptional with respect to $(\eps_{u,i})$.
\end{proof}

\section{A simple $\ell^1$ bound}

This section is devoted to a preliminary bound
for simplified multilinear operators. 
For any tree $T$ and any sequences $y_v\in\ell^1$, define
\begin{equation}
\tilde S_T(y_v)_{v\in\Tterm}(n) = \sum_{\jarrow: j_{v_0}=n}^\star 
\ \prod_{u\in \Tterm} y_u(j_u).
\end{equation}
The notation
$\sum_{\jarrow: j_{v_0}=n}^\star$ indicates that
the sum is taken over all indices $\jarrow\in \integers^T$ that satisfy 
\eqref{cyclicity} as well as $j_{v_0}=n$; 
the restrictions \eqref{exclusion} and \eqref{simplicity} are not imposed here.

\begin{lemma} \label{lemma:Lone}
For any tree $T$ and any sequences 
$\{(y_v): v\in \Tterm\}$
\begin{equation}
\|\tilde S_T(y_v)_{v\in\Tterm}\|_{\ell^1} \le  \prod_{w\in \Tterm}\|y_w\|_{\ell^1},
\end{equation}
with equality when all $y_v(j_v)$ are nonnegative.
\end{lemma}

\begin{proof}
Recall that for some nonnegative integer $k$,
$|T|=3k+1$, $|\Tterm|=2k+1$, and $| \Tnon| = k$.
Consider the set $B\subset T$  whose elements are
$v_0$ together with all $(v,i)$ 
such that $v\in \Tnon$ and $i\in\{1,3\}$.
Thus $|B|=1+2k = |\Tterm|$.
Define
\begin{equation}
k_{v,i} = j_v-j_{(v,i)} \text{ for $v\in \Tnon$ and $i\in\{1,3\}$.}
\end{equation}
Consider the $\integers$-linear mapping $L$ from 
$\integers^{\Tterm}$ to $\integers^B$ defined so that
$L(\jarrow)$ has coordinates $j_{v_0}$ and all $k_{v,i}$.

$j_v$ and $j_{(v,i)}$ are well-defined linear functionals
of $\jarrow\in\integers^{\Tterm}$, because given the quantities
$j_w$ for all $w\in\Tterm$, $j_v$ can be recovered for all
other $v\in T$ via the relations \eqref{cyclicity}, by ascending induction on $v$.
We claim that $L$ is invertible. 
Indeed, from the quantities $j_{v_0}$ and all $j_v-j_{(v,i)}$ with
$v\in \Tnon$ and $i\in\{1,3\}$,
$j_u$ can be recovered for all $u\in T$ by descending induction on $u$,
using again \eqref{cyclicity} at each stage. For instance, at the initial step,
$j_{(v_0,i)} = j_{v_0} + k_{v_0,i}$  for $i=1,3$,
and then $j_{(v_0,2)}$ can be recovered via \eqref{cyclicity}.
Thus $L$ is injective, hence invertible.


By descending induction on nodes it follows in the same way from \eqref{cyclicity}
that $\jarrow=(j_w)_{ w\in\Tterm}$ satisfies a certain linear relation of the 
form 
\begin{equation}j_{v_0} = \sum_{w\in\Tterm} \pm_w j_w\end{equation}
where each coefficient $\pm_w$ equals $\pm 1$. 
By the conclusion of the preceding paragraph, this can be the only relation
to which $(j_w)_{w\in\Tterm}$ is subject; the sum defining 
$\tilde S_T(y_w)_{w\in\Tterm}(j_{v_0})$
is taken over all $\jarrow$ satisfying this relation.
Therefore 
$\sum_{j_{v_0}}\tilde S_T(j_{v_0})$ equals the summation over
all $w\in\Tterm$ and all $j_w\in\integers$, without restriction,
of $\prod_{w\in\Tterm} y_w(j_w)$.
The lemma follows.
\end{proof}

\begin{corollary}
For any ornamented tree,
the sum defining $\tilde\scripts_T(y_v)_{v\in\Tterm}(n)$
converges absolutely for all $n\in\integers$ whenever all $y_v\in\ell^1$, 
and the resulting
sequence satisfies
\begin{equation}
\|\tilde\scripts_T(y_v)_{v\in\Tterm}\|_{\ell^1}
\le 
\prod_{v\in \Tterm} \|y_v\|_{\ell^1}.
\end{equation}
\end{corollary}

There is no bound for $\tilde S_T$ in terms of the quantities
$\|y_w\|_{\ell^p}$ for $p>1$.
It is the additional factors $\langle \rho_u\rangle^{-1}$ 
in the second tree coefficient bound \eqref{maincoefficientbound},
reflecting the dispersive character of the
partial differential equation, which make possible estimates in terms of
weaker $\ell^p$ norms.

\begin{proof}[Proof of Proposition~\ref{prop:formalexpansion}]
Consider any tree $T$ and associated function
\begin{equation} \label{uglybutnecessary}
\int_{\scriptr(T,t)}
\sum_{\jarrow\in\scriptj(T)} 
\prod_{v\in\Tnon}e^{\pm_v i\sigma_v t_v}
\prod_{u\in\Tterm} y_u(t_u,j_u)\,dt_u
\end{equation}
for $0\le t\le\tau$, with $t_{v_0}\equiv t$,
under the assumption that
for each $u\in\Tterm$, $y=y_u$ belongs to
$C^0([0,\tau],\ell^1)$. 
Here for each $t_u$, $y_u(t_u)$
is the sequence whose components are $y_u(t_u,j_u)_{j_u\in\integers}$.

Let each node $v\in\Tterm$ be designated as either finished or unfinished.
Assume that for each finished node, $y_u(t_u,j_u)$ is independent of
$t_u$, while
for each unfinished node,
either the sequence-valued function $t_u\mapsto y_u$ satisfies
the integral equation
\begin{multline}  \label{firstintegraleqnrestated}
y_u(t,n) = y_u(0,n)
-i\omega\int_0^t |y_u(s,n)|^2 y_u(s,n)\,ds
\\
+ i\omega\sum_{j-k+l=n}^* 
\int_0^t y_u(s,j) \bar y_u(s,k) y_u(s,l)e^{i\sigma(j,k,l,n)s}\,ds,
\end{multline}
or its complex conjugate satisfies this same equation.

The $C^0(\ell^1)$ hypothesis guarantees that if we substitute 
the right-hand side of \eqref{firstintegraleqnrestated} for 
$y_u(t_u,j_u)$ in \eqref{uglybutnecessary} for each unfinished node, 
then an absolutely convergent integral and sum are obtained.
Thus we may interchange the outer integral with the sums.
What results is a finite linear combination of 
expressions of the same character as \eqref{uglybutnecessary},
each associated to a larger tree $T^\dagger\supset T$.
At most $3^{|\Tterm|}$ such expressions are obtained, and each is
multiplied by a unimodular numerical coefficient.
Each non-terminal node of $T$ is a non-terminal node of $T^\dagger$,
each finished node of $\Tterm$ remains a terminal node of $T^\dagger$,
and each unfinished node of $\Tterm$ becomes a non-terminal node
of $T^\dagger$, each of whose three children may independently be
either finished or unfinished.
An unfinished node $u$ of $T$ gives rise either to a simple node or a general
node of $T^\dagger$, depending on which of the two trilinear terms
of \eqref{firstintegraleqnrestated} is substituted 
for $y_u$ in \eqref{uglybutnecessary}.

This discussion justifies the formal derivation of the expansion 
in \S\ref{section:ode}. It follows by recursion that for any solution
$u(t,x)$ of the modified Cauchy problem \eqref{modnlsivp}
in $C^0([0,\tau],H^s)$ for sufficiently large $s$,
the associated coefficients $a_n(t) = e^{in^2 t}\widehat{u}(t,n)$
are given by the absolutely convergent infinite power series
\eqref{partofformalsolutionseries}, \eqref{formalsolutionseries} 
for all sufficiently small $t$.
\end{proof}

\section{Tree sum majorants}

\begin{definition}
Let $T$ be an ornamented tree.
The tree sum majorant associated to $T$ is the multilinear operator
\begin{equation}
S_T(y_w)_{w\in\Tterm}(n) = \sum_{\jarrow\in\scriptj(T): j_{v_0}=n}
\ \prod_{u\in \Tnon}
\langle \rho_u(\jarrow)\rangle^{-1} \prod_{w\in \Tterm} y_w(j_w) .
\end{equation}
\end{definition}
Here $t\ge 0$ and $S_T$ is initially defined when all $y_w\in\ell^1$, 
in order to ensure absolute convergence of the sum. 

\begin{definition}
Let $(T,T')$ be a weathered ornamented tree.
The associated tree sum majorant is the multilinear operator
\begin{equation}
S_{(T,T')}(y_w)_{w\in\Tterm}(n) = 
\sum_{\jarrow\in\scriptj(T,T'): j_{v_0}=n}
\ \prod_{u\in \Tnon}
\langle \rho_u(\jarrow)\rangle^{-1} \prod_{w\in \Tterm} y_w(j_w). 
\end{equation}
\end{definition}
Thus
\begin{equation}
S_T = \sum_{T'\subset\Tnon}S_{(T,T')},
\end{equation}
the sum being taken over all subsets $T'\subset  \Tnon$.
The total number of such subsets is $2^{|\Tnon|}\le 2^{|T|}
\le 2^{3|\Tterm|/2} = C^{|\Tterm|}$.

Let $(T,T')$ be a weathered ornamented tree. 
We seek an upper bound
for the associated tree sum operator $S_{(T,T')}$. 
The factors $\langle\rho_v\rangle^{-1}$ 
in the definition of $S_{(T,T')}$ are favorable when $|\rho_v|$ is large;
frozen nodes are those for which $|\rho_v|$ is relatively small,
and hence these require special attention.

Denote by $\Gamma=(\gamma_u)_{ u\in T'}$ any element of $\integers^{T'}$.
Let 
\begin{equation}
\scriptj(T,T',\Gamma)
=\{\jarrow\in\scriptj(T,T'): \rho_u(\jarrow)=\gamma_u
\text{ for all } u\in T' \}.
\end{equation}
$T'$ is the set of all frozen nodes, so by its definition we have
\begin{equation} \label{gammarestriction}
|\gamma_u|=|\rho_u(\jarrow)| 
\le c_0|\sigma_u(\jarrow)|^{1-\delta} \ \forall u\in T'
\end{equation}
with the shorthand notation
$\sigma_u(\jarrow) = \sigma(j_{(u,1)},j_{(u,2)},j_{(u,3)},j_u)$ introduced earlier.
In the remainder of the discussion, we always assume tacitly that
$\Gamma$ satisfies \eqref{gammarestriction}.

This leads to a further decomposition
\begin{multline}
S_{(T,T')}(y_v)_{v\in\Tterm}(n)
= \sum_{\Gamma}
\sum_{\jarrow\in\scriptj(T,T',\Gamma): j_{v_0}=n}\ 
\prod_{u\in\Tnon}\langle\rho_u(\jarrow)\rangle^{-1}
\prod_{w\in \Tterm} y_w(j_w)
\\\
\le 
C^{|T|}
\sum_{\Narrow}
\prod_{v\in T'}2^{-N_v}
\sum_{\Marrow}
\prod_{u\in \Tnon\setminus T'} 2^{-(1-\delta)M_u}
\sum_{\Gamma}
\,\,
\sum_{\jarrow\in\scriptj(T,T',\Gamma): j_{v_0}=n}\ 
\prod_{w\in \Tterm} y_w(j_w)
\end{multline}
where $\Narrow=(N_v)_{v\in T'}$
and $\Marrow=(M_u)_{u\in \Tnon\setminus T'}$.
The notation in the last line means that
the first two sums are taken over all nonnegative integers $N_v,M_u$ as $v$ ranges
over $T'$ and $u$ over $\Tnon\setminus T'$;
the third sum is taken over all $\Gamma$ such that
\begin{equation} \label{gammasizerestriction}
\langle\gamma_v\rangle \in [2^{N_v},2^{1+N_v}) \ \text{ for all $v\in T'$;}
\end{equation}
and the sum with respect to $\jarrow$ is taken over all $\jarrow$ satisfying
the additional restrictions 
\begin{align}
&|\sigma_u(j_{(u,1)},j_{(u,2)},j_{(u,3)},j_u)|\sim 2^{M_u}
\ \text{ for all $u\in \Tnon\setminus T'$}
\label{Marrowrestriction}
\\
&\text{$\rho_v(\jarrow)=\gamma_v$ for all $v\in T'$.}
\label{jarrowrestriction}
\end{align}
Thus 
there is an upper bound $2^{N_v} \le Cc_0|\sigma_v(\jarrow)|^{1-\delta}$   
for all $v\in T'$.

For any $v\in T'$ and any parameter $\gamma_v$,
for any $\jarrow\in\scriptj(T,T',\Gamma)$,
$\sigma_v(\jarrow) =\gamma_v -\sum_{i=1}^3 \eps_{v,i} \rho_{(v,i)}$
where $\rho_{(v,i)}=\rho_{(v,i)}(\jarrow)$ depends only
on $\{j_w-j_{(w,i)}: w<v, i\in\{1,2,3\}\}$. 
Since the quantity $\sigma_v$ on the left-hand side equals
$2(j_v-j_{(v,1)})(j_v-j_{(v,3)})$, for any $\{j_w-j_{(w,l)}: w<v, l\in\{1,2,3\}\}$ 
and any $\gamma_v$
there are at most 
$C_{\delta_1} |\gamma_v-\sum_{i=1}^3 \eps_{v,i} \rho_{(v,i)}|^{\delta_1}$
ordered pairs $\big(j_v-j_{(v,1)},j_v-j_{(v,3)} \big)$
satisfying \eqref{jarrowrestriction}.
Here $\delta_1$ is an arbitrarily small  constant, to be chosen later.

For any frozen node $v\in T'$, $|\gamma_v|$ is small relative
to $\sum_{i=1}^3|\rho_{(v,i)}|^{1-\delta}$, provided that $c_0$ is taken to be small
in the definition of a frozen node.
Therefore we can choose for each combination of parameters
$\Narrow,\Marrow$ a family $\scriptf=\scriptf_{\Narrow,\Marrow}$ 
of vector-valued functions
$F=(f_{v,i}: v\in T', i\in\{1,3\})$
of cardinality at most 
$C_{\delta_1}^{|T|} \prod_{v\in T'} 2^{\max_i N_{(v,i)}\delta_1}$
such that for any $\Gamma$ satisfying \eqref{gammasizerestriction}
and any $\jarrow\in\scriptj(T,T',\Gamma)$,
there exists $F\in\scriptf_{\Narrow,\Marrow}$ 
such that for each $v\in T'$ and each $i\in\{1,3\}$,
\begin{equation} \label{choicefunctionF}
k_{v,i} = j_v-j_{(v,i)} = f_{v,i}(\gamma_v,(k_{w,i}: w<v)).
\end{equation}

Thus for all nonnegative sequences $y_w$ and all $n\in\integers$,
\begin{equation} \label{exhausted1}
|S_{(T,T')}(y_w)_{w\in\Tterm}(n)|
\le 
C^{|T|}
\sum_{\Narrow,\Marrow} 2^{-|\Narrow|}2^{-(1-\delta)|\Marrow|}
\sum_{\Gamma} \sum_{F\in\scriptf_{\Narrow,\Marrow}}
|S_{T,T',\Narrow,\Marrow,\Gamma,F}(y_w)_{w\in\Tterm}(n)|
\end{equation}
where
$|\Narrow|=\sum_v N_v$, $|\Marrow| = \sum_u M_u$, and
\begin{equation} \label{exhausted2}
S_{T,T',\Narrow,\Marrow,\Gamma,F}(y_w)_{w\in\Tterm}(n)
=
\sum_{\jarrow\in\scriptj(T,T',\Gamma): j_{v_0}=n}
\ 
\prod_{w\in \Tterm} y_w(j_w).
\end{equation}
In \eqref{exhausted1},
the second summation is taken over all $\Gamma=(\gamma_u)_{ u\in T'}$
satisfying both \eqref{gammasizerestriction} and \eqref{gammarestriction}.
In \eqref{exhausted2},
the sum is taken over all $\jarrow\in\scriptj(T,T',\Gamma)$
satisfying $j_{v_0}=n$,
\eqref{Marrowrestriction}, \eqref{jarrowrestriction}, and
the additional restriction \eqref{choicefunctionF}. 
We have finally arrived at our basic building blocks,
the multilinear operators $S_{T,T',\Narrow,\Marrow,\Gamma,F}$.

\begin{lemma} \label{lemma:ellpbound}
Let $p\in[1,\infty)$ and $\delta_1>0$.
Then for every exponent $q$ satisfying $q\ge 1$ and $q>p/|\Tterm|$,
there exists $C<\infty$  such that
for every $T,T',\Narrow,\Marrow,\Gamma,F$
and for every sequence $y_v$,
\begin{equation} \label{blockbound}
\| S_{T,T',\Narrow,\Marrow,\Gamma,F}(y_v)_{ v\in \Tterm} \|_{\ell^q}
\le C^{|T|} 
2^{(1+\delta_1)|\Marrow|}
\prod_{v\in\Tterm}\|y_v\|_{\ell^p}.
\end{equation}
\end{lemma}

\begin{proof}
As was shown in the proof of Lemma~\ref{lemma:Lone},
each quantity $j_v$ in the summation defining 
$S_{T,T',\Narrow,\Marrow,\Gamma,F}(y_w)_{ w\in \Tterm}(j_{v_0})$
can be expressed as a function, depending on $\Gamma,F$,
of $j_{v_0}$ together with all $k_{w,i}=j_w-j_{(w,i)}$, where $w$
varies over the set $\Tnon\setminus T'$
of all nodes that are neither frozen nor terminal,
and $i$ varies over $\{1,3\}$.
More precisely, $j_v$ equals $j_{v_0}+g_v$, where $g_v$ is some
function of all these $k_{w,i}$. 

$\prod_{v\in\Tterm}y_v(j_v)$
can thus be rewritten as $\prod_{v\in\Tterm} y_v(j_{v_0}+g_v)$.
If every $k_{w,i}$ is held fixed, then 
as a function of $j_{v_0}$, this product belongs to $\ell^{q}$
for $q=p/|\Tterm|$ with bound $\prod_{v\in\Tterm}\|y_v\|_{\ell^p}$, 
by H\"older's inequality.

The total number of terms 
in the sum defining $S_{T,T',\Narrow,\Marrow,\Gamma,F}$
is the total possible number of vectors $(k_{w,i})$ where
$w$ ranges over $\Tnon\setminus T'$ and $i$ over $\{1,3\}$. 
The number of such pairs
for a given $w$ is $\le C_{\delta_1} 2^{(1+\delta_1)M_w}$,
since $2|k_{w,1}k_{w,3}|=|\sigma_w(\jarrow)|\le 2^{M_w+1}$.
Thus in all there are at most
$C_{\delta_1}^{|T|}2^{(1+\delta_1)|\Marrow|}$ terms.
Minkowski's inequality thus gives the stated bound.
\end{proof}

\eqref{blockbound} is a satisfactory bound, but 
it must be summed over all $F\in\scriptf_{\Narrow,\Marrow}$.
An upper bound for the number 
of such functions $F$ is, roughly speaking,
$C_{\delta_1}^{|T'|}$ 
times the product over all $w\in T'$ of
$\max_i|\rho_{(w,i)}|^{\delta_1}$.
However, this does not quite make sense since $\max_i|\rho_{(w,i)}|$
is a function of $\jarrow$, which we wish to allow to vary while 
$F\in\scriptf_{\Narrow,\Marrow}$ remains fixed.
Thus a correct upper bound is 
\begin{equation}
 |\scriptf_{\Narrow,\Marrow}|
 \le 
 C_{\delta_1}^{|T'|}
 \prod_{v\in T'} 
 2^{\max_i K_{(v,i)}\delta_1}
\end{equation}
where $K_u=N_u$ for $u\in T'$ and $K_u=M_u$ for $u\in\Tnon\setminus T'$,
and the maximum is taken over $i\in\{1,3\}$.

A difficulty now appears.
For each $v\in T'$ we have a compensating
factor of $\langle\gamma_w(\jarrow)\rangle^{-1}\sim 2^{-N_v}$, 
but there is no upper bound whatsoever for the ratio
$\max_i|\rho_{(v,i)}|^{\delta_1}\,/\, 
\langle\gamma_v\rangle$. Thus the factor gained for a given $v\in T'$
cannot compensate for the factor lost for that same node. 
However in aggregate the factors gained compensate for those lost,
as will now be shown.
\begin{lemma} \label{scriptfbound}
For any $\eps>0$ there exists $C_\eps<\infty$ such that uniformly for
all $T,T',\Narrow,\Marrow$,
\begin{equation}
|\scriptf_{\Narrow,\Marrow}|
\le
C_\eps^{|T|}
2^{\eps|\Marrow|}.
\end{equation}
\end{lemma}

\begin{proof}
If the constant $c_0$ in the definition \eqref{frozennode}
of a frozen node is chosen to be sufficiently small, then
any frozen node $u$ has a child $(u,i)$ such that
$|\rho_u|\le \tfrac12|\rho_{(u,i)}|^{1-\delta}$.
Consider any chain $v=u_h\ge u_{h-1}\ge\cdots\ge u_1$
of nodes such that 
$u_{k+1}$ is the parent of $u_k$ for each $1\le k<h$
($u_k$ is called the $(k-1)$-th generation ancestor of $u_1$),
$u_k$ is frozen for all $k>1$, $u_1$ is either not frozen
or is terminal,
and $|\rho_{u_k}| \le \tfrac12|\rho_{u_{k-1}}|^{1-\delta}$.
Then $|\rho_{u_k}|\le 2^{1-k}|\rho_{u_1}|^{(1-\delta)^{k-1}}$;
hence
$2^{K_{u_k}}\le C2^{M_{u_1}(1-\delta)^{k-1}}$.

If $u_1$ is terminal then $\rho_{u_1}=0$ by definition, whence
the inequality $|\rho_{u_k}|\le 2^{1-k}|\rho_{u_1}|^{(1-\delta)^{k-1}}$
forces $\rho_{u_k}=0$ for all $u_k$, as well.
This means that $2^{\max_i K_{(u_k,i)}}\sim 1$. In particular,
this holds for $u_k=v$, so the factor $2^{\max_i K_{(v,i)}}$
will be harmless in our estimates. 
We say that a node $v$ is negligible if there exists such a chain, with $v=u_h$
for some $h\ge 1$.

For each nonnegligible frozen node $v$, choose one such 
chain with $u_h=v$, thus uniquely specifying $h$ and $u_1$ as functions
of $v$;
we then write $u_1=D(v)$.
Given $u_1$ and $h$, there can be at most one $v$ such that $u_1=D(v)$
and $v$ is the $h$-th generation ancestor of $u_1$, simply because 
any node has at most one $h$-th generation ancestor. 
Now taking the product only over nonnegligible nodes $v\in T'$ 
on the left-hand side,
\begin{equation}
\prod_{v\in T'}
2^{ \max_i K_{(v,i)}\delta_1}
\le
\prod_{w\in \Tnon\setminus T'} \prod_{h=1}^\infty
2^{(1-\delta)^{h-1}\delta_1 M_w} 
=
\prod_{w\in \Tnon\setminus T'}
2^{M_w \delta_1/\delta},
\end{equation}
since\footnote{
The exponent $1-\delta<1$ in the definition \eqref{frozennode} of a frozen
node was introduced solely in order to produce a summable series of exponents
$(1-\delta)^{h-1}\delta_1$ in this argument.}
each factor $2^{\max_i K_{(v,i)}\delta_1}$ in the first product is majorized by
$2^{(1-\delta)^{h-1}\delta_1M_w}$ in the second product, where $w=D(v)$
and $v$ is the $h$-th generation ancestor of $w$.
This is not so for negligible nodes, but they contribute at most $C^{|T|}$
to the left-hand side so the conclusion remains valid for the full product.
Thus by choosing $\delta_1$ so that $\delta_1/\delta=\eps$,
since $|\scriptf_{\Narrow,\Marrow}|
\le C_{\delta_1}^{|T|}\prod_{v\in T'}2^{\max_i K_{(v,i)}\delta_1}$,
we obtain $|\scriptf_{\Narrow,\Marrow}|\le C_{\eps}^{|T|}  \prod_{w\in \Tnon\setminus T'}
2^{\eps M_w}=C_{\eps}^{|T|}2^{\eps|\Marrow|}$.
\end{proof}

\begin{proof}[Conclusion of proof of Proposition~\ref{prop:treeoperatorbound}]
Combining the preceding two lemmas gives
\begin{equation}
\sum_{F\in\scriptf_{\Narrow,\Marrow}}\|S_{T,T',\Narrow,\Marrow,\Gamma,F}
(y_v)_{v\in\Tterm}\|_{\ell^q}
\le C_\eps^{|T|} 2^{(1+\eps)|\Marrow|}\prod_{v\in\Tterm}\|y_v\|_{\ell^p}
\end{equation}
for arbitrarily small $\eps>0$,
provided $q\ge\max(1,\frac{p}{|\Tterm|})$.
Since $|\Gamma|\le C^{|T|}2^{|\Narrow|}$, it follows that
\begin{equation}
\sum_{\Gamma}
\sum_{F\in\scriptf_{\Narrow,\Marrow}}\|S_{T,T',\Narrow,\Marrow,\Gamma,F}
(y_v)_{v\in\Tterm}\|_{\ell^q}
\le C_\eps^{|T|} 2^{|\Narrow|}2^{(1+\eps)|\Marrow|}\prod_{v\in\Tterm}\|y_v\|_{\ell^p}.
\end{equation}
On the other hand, Lemma~\ref{lemma:Lone}
gives a uniform $\ell^1$ norm bound of $C^{|T|}
\prod_{v\in\Tterm}\|y_v\|_{\ell^1}$  for the summation over all $\jarrow$.
Thus if $q>\frac{p}{|\Tterm|}$ and $q\ge 1$ we may interpolate to find
that there exists $\eta>0$ depending on $q-\frac{p}{|\Tterm|}$
but not on $\delta$
such that
\begin{equation}
\sum_{\Gamma}
\sum_{F\in\scriptf_{\Narrow,\Marrow}}\|S_{T,T',\Narrow,\Marrow,\Gamma,F}
(y_v)_{v\in\Tterm}\|_{\ell^q}
\le C_\eta^{|T|} 2^{(1-\eta)|\Narrow|+(1-\eta)|\Marrow|}
\prod_{v\in\Tterm}\|y_v\|_{\ell^p}.
\end{equation}
Taking into account the factors $2^{-|\Narrow|}2^{-(1-\delta)|\Marrow|}$
in \eqref{exhausted1}, 
summing over $\Narrow,\Marrow$
as well as over all subsets $T'\subset\Tnon$
yields a convergent series and
completes the proof of Proposition~\ref{prop:treeoperatorbound}.
\end{proof}

\section{Loose ends}
We may reinterpret the sum of our power series
\eqref{partofformalsolutionseries},\eqref{formalsolutionseries}
as a function via the relation $\widehat{u}(t,n) = e^{in^2 t} a_n(t)$
with $a(0)$ defined by $\widehat{u_0}(n) = a_n(0)$,
and will do so consistently without further comment,
abusing notation mildly by writing $u(t,x) = S(t)u_0(x)$.

\begin{lemma}
Let $p\in[1,\infty)$. For any $R>0$ there exists $\tau>0$ such that
for any $u_0\in \scripth^p$ with norm $\le R$,
the element $u(t,x)\in C^0([0,\tau],\scripth^p)$ defined
by \eqref{partofformalsolutionseries},\eqref{formalsolutionseries}
is a limit, in $C^0([0,\tau],\scripth^p)$ norm,
of smooth solutions of \eqref{modnlsivp}. 
\end{lemma}

\begin{proof}
All of our estimates apply also
in the spaces $\scripth^{s,p}$
defined by the condition that 
$(\langle n\rangle^s\widehat{f}(n))_{n\in\integers}\in\ell^p$,
provided that $1\le p<\infty$ and $s>0$.
This follows from the proof given for $s=0$ above, for
the effect of working in $\scripth^{s,p}$ is to introduce a factor of
$\prod_{v\in\Tnon} 
\frac{\langle j_v\rangle^{s}}{\prod_{i=1}^3 
\langle j_{(v,i)} \rangle^s}$
in the definition of the tree operator. 
The relation \eqref{cyclicity} ensures that $\max_i |j_{(v,i)}|\ge\tfrac13|j_v|$,
whence $\frac{\langle j_v\rangle^{s}}{\prod_{i=1}^3 \langle j_{(v,i)} \rangle^s}
\lesssim 1$, so the estimates for $s=0$ apply directly to all $s>0$.
More generally, if $\scripth^{s,p}$ is equipped with the norm
\begin{equation*}
\|f\|_{\scripth^{s,p}_\eps}
= \|(1+|\eps\cdot|^{2s})^{1/2}\widehat{f}(\cdot)\|_{\ell^p}
\end{equation*}
then all estimates hold uniformly in $\eps\in[0,1]$
and $s\ge 0$.

Given $s$ and any initial datum $u_0$ satisfying $\|u_0\|_{\scripth^p}\le R$
with the additional property that $\widehat{u_0}(n)=0$ for all 
$|n|>N$,
we may choose $\eps>0$ so that 
$\|u_0\|_{\scripth^{s,p}_\eps}\le 2R$;
$\eps$ depends on $N$ but not on $R$.
Thus the infinite series converges absolutely and uniformly
in $C^0([0,\tau],H^{s-\tfrac12+\tfrac1p})$ 
if $p\ge 2$ and in $C^0([0,\tau],H^{s})$ if $p\le 2$,
where $\tau$ depends only on $R$, not on $s$.
By Lemma~\ref{automaticsolution}, the series sums to a solution
of \eqref{modnlsivp} in the sense \eqref{ufirstintegraleqn}; but since
the sum is very smooth as a function of $x$
(that is, its Fourier coefficients decay rapidly)
this implies that it is a solution in the classical sense.
Given an arbitrary $u_0$ satisfying $\|u_0\|_{\scripth^p}\le R$,
we can thus approximate it by such special initial data
to conclude that $S(t)u_0$
is indeed a limit, in $C^0([0,\tau],\scripth^p)$, of smooth solutions. 
\end{proof}
 
\begin{proof}[Proof of Proposition~\ref{prop:nonlineardistribution}]
Let $u_0\in\scripth^p$ be given,
let $u(t,x)=S(t)(u_0)\in C^0([0,\tau],\scripth^p)$.
We aim to prove that the nonlinear term
$\omega|u|^2 u$ has an intrinsic meaning as 
$\lim_{N\to\infty} \omega |T_N u|^2 T_Nu$
in the sense of distributions in $(0,\tau)\times\torus$.
Forming $T_NS(t)(u_0)$ is of course not the same thing as forming
$S(t)(T_Nu_0)$.

Define $a_n(t) = e^{in^2 t} \widehat{u}(t,n)$. 
Denote also by $T_N$ the operator that maps a sequence-valued function $(b_n(t))$
to $(T_Nb_n(t))$ where $T_N b_n=b_n$ if $|n|\le N$, and $=0$ otherwise.
It suffices to prove that 
\begin{equation} 
\int_0^t 
\sum_{j-k+l=n}^* 
T_Na_{j}(s) \overline{T_Na_{k}(s)} T_Na_{l}(s)e^{i\sigma(j,k,l,n)s}\,ds
-\int_0^t |T_Na_{n}(s)|^2 T_Na_{n}(s)\,ds
\end{equation}
converges
in $\ell^p$ norm as $N\to\infty$, uniformly for all $t\in[0,\tau]$,
to
\begin{equation*} 
\sum_{j-k+l=n}^* 
\int_0^t 
a_{j}(s) \overline{a_{k}(s)} a_{l}(s)e^{i\sigma(j,k,l,n)s}\,ds
-\int_0^t |a_{n}(s)|^2 a_{n}(s)\,ds.
\end{equation*}
Convergence in the distribution
sense follows easily from this by expressing any sufficiently smooth function
of the time $t$ as a superposition of characteristic functions of intervals $[0,t]$.

Now in the term 
$ \int_0^t 
\sum_{j-k+l=n}^* 
T_Na_{j}(s) 
\overline{T_Na_{k}(s)} 
T_Na_{l}(s)e^{i\sigma(j,k,l,n)s}\,ds$,
the integral may be interchanged with the sum since the truncation operators
restrict the summation to finitely many terms.
Expanding $a_j,a_k,a_l$ out as infinite series of tree operators
applied to $a(0)$, we obtain finally an infinite series 
of the general form $\sum_{k=1}^\infty B_k(t)(a(0),\cdots,a(0))$
where $B_k(t)$ is a finite linear combination of $O(C^k)$ tree sum operators,
with coefficients $O(C^k)$,
applied to $a(0)$ just as before, 
with the sole change that the extra restriction
$|j_{(v_0,i)}|\le N$ for $i\in\{1,2,3\}$
for indices corresponding to children of the root node
is placed on $\jarrow$ in the summation defining $\scripts_T$
for each tree $T$.

Since we have shown that all bounds hold for the sums of the absolute
values of the terms in the tree sum, 
it follows immediately that this trilinear term converges as $N\to\infty$. 
Convergence for the other nonlinear term is of course trivial.
Likewise it is trivial that $(T_N u)_t\to u_t$ and $(T_N u)_{xx}\to u_{xx}$,
by linearity.
\end{proof}
This reasoning shows that the limit of each term equals the sum of a convergent
power series, taking values in $C^0([0,\tau],\scripth^p)$, in $u_0$.

Given $R>0$, there exists $\tau>0$ for which we have shown that
for any $a(0)\in\ell^p$ satisfying $\|a(0)\|_{\scripth^p}\le R$,
our power series expansion defines $a(t)\in C^0([0,\tau],\ell^p)$,
as an $\ell^p$-valued analytic function of $a(0)$.
Moreover for any $t\in[0,\tau]$, 
both cubic terms in the integral equation \eqref{firstintegraleqn}
are well-defined as limits obtained by replacing $a(s)$ by $T_N a(s)$,
evaluating the resulting cubic expressions, and passing to the limit
$N\to\infty$.

\begin{lemma}
Whenever $\|a(0)\|_{\ell^p}\le R$,
the function $a(t)\in C^0([0,\tau],\ell^p)$
defined as the sum of the power series expansion \eqref{partofformalsolutionseries}
satisfies the integral equation \eqref{ufirstintegraleqn} when 
the nonlinear terms in \eqref{firstintegraleqn} are defined by the limiting procedure
described in the preceding paragraph.
\end{lemma}

\begin{proof}
This follows by combining Lemma~\ref{automaticsolution}
with the result just proved.
\end{proof}

\begin{proof}[Proof of Proposition~\ref{prop:weaksmoothing}]
Let $u_0\in\scripth^p$.
If $u=Su_0$, and if $v$ is the solution of the Cauchy problem \eqref{modnlsivp} for
the modified linear Schr\"odinger equation with initial datum $u_0$,
then $u_0-v$ is expressed as  $\sum_{k=1}^\infty B_k(t)(u_0,\cdots,u_0)$
where the $n$-th Fourier coefficient of $B_k(t)(u_0,\cdots)(t)$ equals
$e^{-in^2 t} A_k(t)(a(0),\cdots)$ with $a_n(0) = \widehat{u_0}(n)$.
According to Proposition~\ref{prop:treeoperatorbound},
$\|A_k(t)(a(0),\cdots)\|_{\ell^q} = O(t^{k\eps}\|a(0)\|_{\ell^p}^k)$
whenever $q>\tfrac{p}3$ and $q\ge 1$. 
Summing over $k$ yields the conclusion.
\end{proof}

\end{document}